\documentclass[11pt]{article}
\usepackage{epic,latexsym,amssymb,xcolor}
\usepackage{color}
\usepackage{tikz}
\usepackage{amsfonts,epsf,amsmath,leftidx}
\usepackage{float}
\usepackage{pgfplots}
\pgfplotsset{compat=1.15}
\usepackage{mathrsfs}
\usetikzlibrary{arrows}

\usepackage[bottom=1.5in, top=1.2in, left=1.5in, right=1.5in]{geometry}

\newtheorem{thm}{Theorem}
\newtheorem{lem}{Lemma}

\newtheorem{cor}{Corollary}
\newtheorem{ob}{Observation}

\newtheorem{problem}{Problem}

\newcommand{\qed}{$\Box$}

\newcommand{\proof}{\noindent\textbf{Proof. }}

\let\oldenumerate\enumerate
\renewcommand{\enumerate}{
  \oldenumerate
  \setlength{\itemsep}{0pt}
  \setlength{\parskip}{0pt}
  \setlength{\parsep}{0pt}
}

\begin{document}

\title{Counting the Numbers of Paths of all Lengths in Dendrimers and its Applications}

\author{Hafsah Tabassum$^{1, 2}$\thanks{Research supported by Petchra Pra Jom Klao Ph.D. Research Scholarship from King Mongkut's University of Technology Thonburi~(1281/2021).} \,,  \,Syed Ahtsham Ul Haq Bokhary$^{3}$,\\
Thiradet Jiarasuksakun$^{1, 2}$, and\, Pawaton Kaemawichanurat$^{1, 2}$\thanks{}
\\ \\

$^{1}$Mathematics and Statistics with Applications (MaSA), \\
$^{2}$Department of Mathematics,\\
King Mongkut's University of Technology Thonburi, \\
Bangkok, Thailand \\
$^{3}$Bahauddin Zakariya University, Multan, Pakistan.\\
\small \tt Email: hafsahtabassum@yahoo.com, sihtsham@gmail.com,\\ thiradet.jia@kmutt.ac.th, pawaton.kae@kmutt.ac.th}

\date{}
\maketitle

\begin{abstract}
For positive integers $n$ and $k$, the dendrimer $T_{n, k}$ is defined as the rooted tree of radius $n$ whose all vertices at distance less than $n$ from the root have degree $k$. The dendrimers are higly branched organic macromolecules having repeated iterations of
branched units that surroundes the central core. Dendrimers are used in a variety of fields including chemistry, nanotechnology, biology. In this paper, for any positive integer $\ell$, we count the number of paths of length $\ell$ of $T_{n, k}$. As a consequence of our main results, we obtain the average distance of $T_{n, k}$ which we can establish an alternate proof for the Wiener index of $T_{n, k}$. Further, we generalize the concept of medium domination, introduced by Varg\"{o}r and D\"{u}ndar in 2011, of $T_{n, k}$.
\end{abstract}


\section{Introduction}
The set of vertices in a graph $G = (V(G), E(G))$ is $V(G)$ while the set of edges is denoted by $E(G)$. All graphs in this paper are finite and simple, with no loops or multiple edges. The the set $\{u : uv \in E(G)\}$ is the \emph{neighbor set} $N_{G}(v)$ of a vertex $v$ in $G$. The degree $deg_{G}(v)$ of a vertex $v$ in $G$ is given by $|N_{G}(v)|$. If the subgraph of $G$ induced by $S$ has no edges, then the vertex subset $S$ of $V(G)$ is \emph{independent}. The maximum cardinality of an independent set is given by the \emph{independence number} of $G$ which denoted by $\alpha(G)$. If every vertex of $G$ has degree $k$ then the graph $G$ is $k$\emph{-regular}. For $u, v \in V(G)$, the length of a shortest path from $u$ to $v$ is the \emph{distance} $d_{G}(u, v)$ between $u$ and $v$ in $G$. The maximum distance between all pairs of vertices of $G$ is \emph{diameter} $diam(G)$. The total of the distance between each pair of vertices of $G$ divided by the number of pairs of vertices is the \emph{average distance} $\mu(G)$ of $G$. That is:

$$\mu(G) = \frac{\displaystyle{\sum_{u, v \in V(G)}d_{G}(u, v)}}{\displaystyle{\binom{|V(G)|}{2}}}.$$

A \emph{tree} is a graph with no subgraphs that are cycles. A \emph{leaf}, also known as \emph{pendant vertex}, is a vertex with degree one. A leaf's incident edge is the \emph{pendant edge} while a leaf's neighbouring vertex is called a \emph{support vertex}. A \emph{rooted tree} $T$ is a tree whose one vertex  identified as the root $r$. Furthermore, if $d_{T}(r, v) = i$  a vertex $v$ of $T$ is at \emph{level} $i$  and $T$ has $n$-\emph{level} if the greatest level of all vertices of $T$ is $n$. A \emph{balanced tree} is a rooted tree with equal number of vertices at the same level having the same degree. A \emph{dendrimer} $T_{n, k}$ is an $n$-level balanced tree with the degree $k$ for all non-leaf vertices.A dendrimer is a molecule with a well-defined chemical structure that is synthesised chemically. Dendrimers have three key main components: one is the core, and it's the most fundamental aspect in dendrimer development, then branches that are added at each step sequentially to produce a structure like tree, the last component is end groups. Dendrimers are hyperbranched macromolecules that have a wide range of applications in domains like supramolecular, drug development, and nanotechnology. Some graph constants such as dominating numbers and some other tpes of dominating number are used to describe a  range of physical characteristics, including physicochemical characteristics, thermodynamic characters, chemical and biological actions, and so on.
In 1978, Fritz Vogtle's was   the first to bring these nanomolecules to researcher's attention \cite{BEV}.
In \cite{SMS} the topological indices of well-known dendrimers were introduced. Researchers have discovered topological indices for many chemical structures such as dendrimers, trees and other graphs, inspired by the chemical relevance of topological  indices of molecular networks. In context of spectral graph theory, the sum of the absolute values of the eigenvalues of the graph $G$ is known as the \emph{energy of graph} $G$ which can be given as $$E(G)= \sum_{i=1}^{n}|\lambda_{i}|$$ where $\lambda_{i}, i=1,2, \dots, n$ are the eigenvalues of $G$.
In \cite{BH}, the eigenvalues of the cayley tree dendrimers are obtained with the help of the characteristic polynomials. The reduction formulas for calculating the characteristic polynomials of d(2, k) and d(3, k) are constructed. Also, The energy of the above mentioned dendrimers is calculated. For more studies in dendrimers see \cite{CDSS,SIA} for example.
\vskip 5 pt

\indent The followings are examples of $T_{n,, k}$ when $n$ or $k$ is small. By the definition of dendrimers, we have that $n \geq 1$ and $k \geq 2$.
\vskip 5 pt

\indent When $n=1$, $T_{1,k}$ can be constructed by introducing $k$ vertices and joining each of them to the central vertex through an edge. Thus  $T_{1, k}$ is a \emph{star} with $n + 1$ vertices.
\vskip 5 pt

\indent When $k=2$, $T_{n,2}$ is a \emph{path} consisting of $2n + 1$ vertices. 
\vskip 5 pt

\indent When $n=2$, it can be observed that $T_{2,k}$ can be constructed from $T_{1,k}$ by introducing $k-1$ vertices at each leaf of $T_{1,k}$ and then joining them to that leaf. Hence, when $n,k \geq 2$, $T_{n,k}$ is constructed from $T_{n-1,k}$ by introducing further $k-1$  vertices at each leaf vertex of $T_{n-1,k}$ and then joining them to that leaf vertex. Namely, the procedure of construction  of $T_{n,k }$ for $(n \geq 2, k \geq 0)$ consists of $n$ iterations from the graph that has exactly one vertex. 
\vskip 30 pt


\begin{figure}[H]
\centering
\definecolor{ududff}{rgb}{0.30196078431372547,0.30196078431372547,1}
\resizebox{0.55\textwidth}{!}{%

\begin{tikzpicture}[line cap=round,line join=round,>=triangle 45,x=1cm,y=1cm]
\draw [line width=0.4pt] (5,4)-- (5,5);
\draw [line width=0.4pt] (5,4)-- (8.17797113184519,4.008784842075203);
\draw [line width=0.4pt] (5,4)-- (1.7926997594517675,3.98661440518);
\draw [line width=0.4pt] (5,4)-- (5.008840390601431,0.82756502783772);
\draw [line width=0.4pt] (5,5)-- (5.010634007491878,7.2003001124159525);
\draw [line width=0.4pt] (6,4)-- (7,5);
\draw [line width=0.4pt] (6,4)-- (7,3);
\draw [line width=0.4pt] (5,3)-- (6,2);
\draw [line width=0.4pt] (5,3)-- (4,2);
\draw [line width=0.4pt] (4,4)-- (3,3);
\draw [line width=0.4pt] (4,4)-- (3,5);
\draw [line width=0.4pt] (5,5)-- (4,6);
\draw [line width=0.4pt] (5,5)-- (6,6);
\draw [line width=0.4pt] (5,6)-- (5.59938,7.00798);
\draw [line width=0.4pt] (5,6)-- (4.3922,6.99603);
\draw [line width=0.4pt] (6,6)-- (6.019389382707226,6.801191539238361);
\draw [line width=0.4pt] (6,6)-- (6.614108955390044,6.5749692336716015);
\draw [line width=0.4pt] (6,6)-- (6.803639066085898,5.995745918471074);
\draw [line width=0.4pt] (7,5)-- (6.994402502583413,5.804982481973559);
\draw [line width=0.4pt] (7,5)-- (7.60684,5.58945);
\draw [line width=0.4pt] (7,5)-- (7.7998481233506975,5.020732798594885);
\draw [line width=0.4pt] (7,4)-- (8.00198,4.59597);
\draw [line width=0.4pt] (7,4)-- (8.00198,3.39929);
\draw [line width=0.4pt] (7,3)-- (7.810446092045004,3.0071187466766687);
\draw [line width=0.4pt] (7,3)-- (7.61813,2.39453);
\draw [line width=0.4pt] (7,3)-- (7.005000471277721,2.191075157215076);
\draw [line width=0.4pt] (6,2)-- (6.814237034780206,2.0003117207175607);
\draw [line width=0.4pt] (6,2)-- (6.59079,1.40105);
\draw [line width=0.4pt] (6,2)-- (6.019389382707227,1.1842681312559677);
\draw [line width=0.4pt] (5,2)-- (5.60027,1.00245);
\draw [line width=0.4pt] (5,2)-- (4.39551,0.98757);
\draw [line width=0.4pt] (4,2)-- (3.9951773620947084,1.1948660999502743);
\draw [line width=0.4pt] (4,2)-- (3.41385,1.40403);
\draw [line width=0.4pt] (4,2)-- (3.179133772633117,1.9791157833289479);
\draw [line width=0.4pt] (3,3)-- (2.988370336135602,2.2016731259093825);
\draw [line width=0.4pt] (3,3)-- (2.40921,2.39995);
\draw [line width=0.4pt] (3,3)-- (2.20412065275693,2.9859228092880556);
\draw [line width=0.4pt] (3,4)-- (2.00357,3.39561);
\draw [line width=0.4pt] (3,4)-- (2.00357,4.58796);
\draw [line width=0.4pt] (3,5)-- (2.20412065275693,5.010134829900578);
\draw [line width=0.4pt] (3,5)-- (2.9989683048299085,5.8155804506678646);
\draw [line width=0.4pt] (3,5)-- (2.40921,5.59592);
\draw [line width=0.4pt] (4,6)-- (3.20032971002173,5.995745918471073);
\draw [line width=0.4pt] (4,6)-- (3.40017,6.58965);
\draw [line width=0.4pt] (4,6)-- (3.9951773620947084,6.80119153923836);
\begin{scriptsize}
\draw [fill=black] (5,4) circle (1.5pt);
\draw [fill=black] (4,4) circle (1.5pt);
\draw [fill=black] (5,5) circle (1.5pt);
\draw [fill=black] (6,4) circle (1.5pt);
\draw [fill=black] (5,3) circle (1.5pt);
\draw [fill=black] (5,6) circle (1.5pt);
\draw [fill=black] (6,6) circle (1.5pt);
\draw [fill=black] (7,5) circle (1.5pt);
\draw [fill=black] (7,4) circle (1.5pt);
\draw [fill=black] (7,3) circle (1.5pt);
\draw [fill=black] (6,2) circle (1.5pt);
\draw [fill=black] (5,2) circle (1.5pt);
\draw [fill=black] (4,2) circle (1.5pt);
\draw [fill=black] (3,3) circle (1.5pt);
\draw [fill=black] (3,4) circle (1.5pt);
\draw [fill=black] (3,5) circle (1.5pt);
\draw [fill=black] (4,6) circle (1.5pt);
\draw [fill=black] (8.17797113184519,4.008784842075203) circle (1.5pt);
\draw [fill=black] (8.00198,4.59597) circle (1.5pt);
\draw [fill=black] (8.00198,3.39929) circle (1.5pt);
\draw [fill=black] (7.60684,5.58945) circle (1.5pt);
\draw [fill=black] (7.7998481233506975,5.020732798594885) circle (1.5pt);
\draw [fill=black] (6.994402502583413,5.804982481973559) circle (1.5pt);
\draw [fill=black] (7.810446092045004,3.0071187466766687) circle (1.5pt);
\draw [fill=black] (7.005000471277721,2.191075157215076) circle (1.5pt);
\draw [fill=black] (7.61813,2.39453) circle (1.5pt);
\draw [fill=black] (6.814237034780206,2.0003117207175607) circle (1.5pt);
\draw [fill=black] (6.019389382707227,1.1842681312559677) circle (1.5pt);
\draw [fill=black] (6.59079,1.40105) circle (1.5pt);
\draw [fill=black] (5.008840390601431,0.82756502783772) circle (1.5pt);
\draw [fill=black] (5.60027,1.00245) circle (1.5pt);
\draw [fill=black] (4.39551,0.98757) circle (1.5pt);
\draw [fill=black] (3.9951773620947084,1.1948660999502743) circle (1.5pt);
\draw [fill=black] (3.179133772633117,1.9791157833289479) circle (1.5pt);
\draw [fill=black] (3.41385,1.40403) circle (1.5pt);
\draw [fill=black] (2.988370336135602,2.2016731259093825) circle (1.5pt);
\draw [fill=black] (2.20412065275693,2.9859228092880556) circle (1.5pt);
\draw [fill=black] (2.40921,2.39995) circle (1.5pt);
\draw [fill=black] (1.7926997594517675,3.98661440518) circle (1.5pt);
\draw [fill=black] (2.00357,3.39561) circle (1.5pt);
\draw [fill=black] (2.00357,4.58796) circle (1.5pt);
\draw [fill=black] (2.9989683048299085,5.8155804506678646) circle (1.5pt);
\draw [fill=black] (2.20412065275693,5.010134829900578) circle (1.5pt);
\draw [fill=black] (2.40921,5.59592) circle (1.5pt);
\draw [fill=black] (3.20032971002173,5.995745918471073) circle (1.5pt);
\draw [fill=black] (3.9951773620947084,6.80119153923836) circle (1.5pt);
\draw [fill=black] (3.40017,6.58965) circle (1.5pt);
\draw [fill=black] (5.010634007491878,7.2003001124159525) circle (1.5pt);
\draw [fill=black] (4.3922,6.99603) circle (1.5pt);
\draw [fill=black] (5.59938,7.00798) circle (1.5pt);
\draw [fill=black] (6.019389382707226,6.801191539238361) circle (1.5pt);
\draw [fill=black] (6.803639066085898,5.995745918471074) circle (1.5pt);
\draw [fill=black] (6.614108955390044,6.5749692336716015) circle (1.5pt);
\end{scriptsize}
\label{t34}
\end{tikzpicture}

 }%
\vskip 1 cm
\caption{The dendrimers $T_{3, 4}$.}
\end{figure}
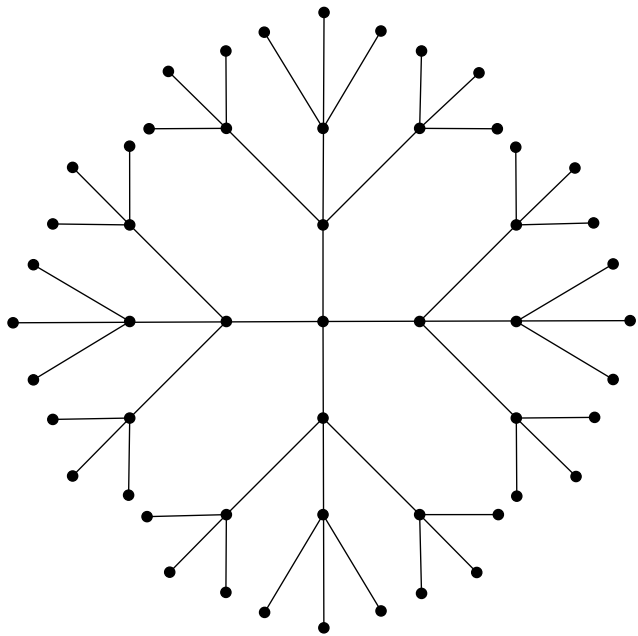

\vskip 5 pt

\indent For a graph $G$, the sum of the distance between any pair of vertices of $G$ is known as the \emph{Weiner index} $W(G)$ of $G$. That is:

$$W(G) = \displaystyle{\sum_{\{u, v\} \subseteq V(G)}d_{G}(u, v)}.$$

\noindent In the quantitative structure-property
relationships (QSPR) \cite{S,N,R}, 
the Wiener index was the first and most well researched topological index. Since then dozens of new indices have been developed to link topological indices with various physical features. The boiling temperatures of alkane molecules are closely associated with the Wiener index number, according to Wiener. Later study on quantitative structure activity linkages revealed it is  connected with some other factors also such as the critical point parameters, density, surface tension, viscosity of the liquid phase and the molecule's van der Waals surface area. It was originally called the path number since it was defined as the total of the lengths  between any two carbon atoms in an alkane in terms of carbon-carbon bonds \cite{WH}. Wiener's works did not make use of graph theory, and the path number was only used in acyclic systems. Hosoya, \cite{HH} in 1971, was the first to define the \emph{Wiener index} within the context of chemical graph theory. The Weiner index appears to have been studied for the first time in the mathematical literature in 1976 \cite{E}. This index has also been referred to by the terms ``graph distance" \cite{E} and ``transmission" \cite{M,P} uses the Laplacian matrix to introduce a graph-theoretical new definition of the Wiener index for trees. Andrey A. Dobryin et. al. established the Wiener index for Cayley tree dendrimer in \cite{DRE} which states that for every $k \geq 3$, the Wiener index of $T_{n,k}$ is 
\begin{align}\label{weiner1}
W(T_{n,k})=\frac{1}{(k-2)^{3}}[(k-1)^{2n}[nk^{3}-2(n+1)k^{2}+k]+2k^{2}(k-1)^{n}-k].
\end{align}
\noindent For more studies of Weiner index see \cite{KST,LL} for example.
\vskip 5 pt

\indent A graph $G$ fulfilling specific constraints can efficiently simulates numerous scenarios in communication, facility locating, cryptology and other fields. Due to cost constraint, it is frequently sought to have a spanning tree of $G$ that is optimal with respect to one or more attributes. One of these attributes is usually the average distance between vertices. The average distance in graph $G$ is defined as follows:
\begin{equation*}
    \mu(G)=\frac{\sum_{\{u, v\} \subseteq V(G)}d_{G}(u, v)}{\binom{|V(G)|}{2}}
\end{equation*}
\noindent The study of the average distance of graphs was initiated by Plesnik in his classical result in \cite{P}. The average distance is an important tool to analyse entire structure of the graph. The parameter globally presents expected number of edges that an object needs to travel between nodes~(vertices) of networks. This reflects data transmission efficiency of communication networks as well as capability to deliver objects of transportation networks. Hence, the average distance has been continuously studied in both theoretical, algorithm and application areas. For example of the studies of average distance of graphs, Fajtlowicz and Waller \cite{FW} established the inequality between the average distance and the independence number in their classical paper since 1986 that $\alpha(G) \geq \mu(G) - 1$ for every connected graph $G$. Chung \cite{C} improved this bound to bt $\alpha(G) \geq \mu(G)$ and further characterized that the equality holds if and only if $G$ is a complete graph. For more studies of the average distance of graphs see \cite{DMS,J,TX} for example.
\vskip 5 pt

\indent Domination in graph has been extensively researched and utilised in a variety of fields. Varg\"{o}r and D\"{u}ndar~\cite{DP}, established the idea of ``the medium domination number" which is defined as the total number of vertices that dominate every pair of vertices with the  average value of it. In the same way each vertex in a graph may protect every vertex in its immediate surroundings and in domination every vertex in neighborhood must be secured. In any connected simple graph $G$ having order $n$, the medium domination number of $G$ is defined as $\gamma_m (G)= \frac {TDT(G)}{{n\choose 2}}$. The medium domination number of Jahangir graph was determind by Ramachandran and Parvathi~\cite{MN}. Mirajkar et. al. found the medium domination number of few poly silicates in \cite{MK}. Mahadevan et.al. proposed the concept of the Extended Medium Domination number of a graph in \cite{GS}. The total number of vertices that dominate each pair of vertices 

$$ETDV (G)=\sum edom (u,v)$$

\noindent for any $u,v \in V(G)$. The extended medium domination number of a graph $G$ is defined as $EMD(G)= \frac{EDT(G)}{{n\choose 2}}$. $ETDV (G)$ 
is the sum of number of path of length one,two and three.
Motivated by the above G. Mahadevan et.al. introduced the idea of Double
Twin domination of a graph \cite{GS} The total number of vertices that dominate every pair of vertices 

$$SDTwin (G)= \sum DTwin (u,v)$$ 

\noindent for $u,v \in V(G)$. In any simple graph $G$ with $n$ number of vertices, the double Twin domination number of $G$ can be given as $DTD(G) = \frac{SDTwin (G)}{{n\choose 2}}$ where $SDTwin (G)$  is the total of number of path having length one, two, three and four. This number was discovered by  Mahadevan and Vijayalakshmi in
\cite{GV3}, for variety of common classes of graphs like path, cycle, wheel graph, complete graph, star graph, Cartesian product of path and Corona product of
path.
\vskip 5pt

\indent From the above discussion, it can be showed that the average distance and the medium domination number of dendrimers can be found if we know the number of paths of all lengths. Thus, the problems that arises is:

\begin{problem}\label{p1}
For non-negative integers $n, k$ and $\ell$, how many paths of length $\ell$ does a dendrimer $T_{n, k}$ have? 
\end{problem}

\noindent Surprisingly, to the best of our knowledge, Problem \ref{p1} has not been answered.
\vskip 5 pt

\indent In this paper, we solve Problem \ref{p1} by establishing the exact and recursive formulas to count the number of paths of length $\ell$ of $T_{n, k}$ for all $1 \leq \ell \leq 2n$. As a consequence, we easily obtain average distance of $T_{n, k}$. Further, we generalize the concept of medium domination to $\ell$-medium domination in graphs.
\vskip 5 pt

\section{Main Results and Applications}
In this section, we state our main results of this paper as well as their applications in Subsections \ref{ap1} and \ref{ap2} while most of the proofs are given in Section \ref{proofs}. First, for a graph $G$, we let
\vskip 5 pt

\indent $n_{\ell}(G)$: the number of paths of length $\ell$ of $G$.
\vskip 5 pt

\noindent The first main result is the formula of $n_{\ell}(T_{n, k})$ for all possible values of $\ell$. Recall that when $k = 2$, the dendrimer $T_{n, 2}$ is a path of length $2n$. Thus, we let $x_{1}, ..., x_{2n + 1}$ be $T_{n, 2}$. Clearly, for a positive integer $1 \leq \ell \leq 2n$, all the paths of length $\ell$ are $x_{i}, x_{i + 1}, ..., x_{i + \ell}$ for all $1 \leq i \leq 2n + 1 - \ell$. Hence, we obtain the following observation.
\vskip 5 pt

\begin{ob}\label{thm1}
Let $T_{n, k}$ be the dendrimer. If $k = 2$, then 
\begin{equation*}
n_{\ell}(T_{n, 2}) = 2n + 1 - \ell.
\end{equation*}
\end{ob}
\vskip 5 pt

\noindent Thus, throughout of this paper, we may assume that $k \geq 3$. Further, for a tree $T$, we let
\vskip 5 pt

\indent $n^{1}_{\ell}(T):$ the number of paths of length $\ell$ of $T$ with exactly one end vertex is a leaf of $T$.
\vskip 5 pt

\indent $n^{2}_{\ell}(T):$ the number of paths of length $\ell$ of $T$ whose both end vertices are leaves of $T$.
\vskip 5 pt

\noindent Our main results in this subsection are Theorem \ref{mainthm}, Corollaries \ref{cor1} and \ref{cor2}. As informed earlier, the proofs are given in Section \ref{proofs}.
\vskip 5 pt

\begin{thm}\label{mainthm}
Let $T_{n, k}$ be the dendrimer and $k \geq 3$. 
If $\ell$ is even number, then 
\begin{equation*}
n_{\ell}(T_{n, k}) = (k - 1)n^{1}_{\ell - 1}(T_{n - 1, k}) + (k - 1)^{2}n^{2}_{\ell - 2}(T_{n - 1, k}) + n_{\ell}(T_{n - 1, k}).
\end{equation*}
If $\ell$ is odd number, then 
\begin{equation*}
n_{\ell}(T_{n, k}) = (k - 1)n^{1}_{\ell - 1}(T_{n - 1, k}) +2 (k - 1)n^{2}_{\ell - 1}(T_{n - 1, k}) + n_{\ell}(T_{n - 1, k}).
\end{equation*}
\end{thm}
\vskip 5 pt

\noindent By Theorem \ref{mainthm}, we obtain the following corollaries. It is worth noting that Corollary \ref{cor2} is a combinatorial identity which is obtained by the counting two way principle.
\vskip 5 pt

\begin{cor}\label{cor1}
Let $T_{n, k}$ be the dendrimer and $k \geq 3$. 
Then
\begin{equation*}
	n_{\ell}(T_{n, k}) = \left\{
	\begin{array}{rl}
		\frac{k(k-1)^{\ell-1}}{2} [\frac{k(k-1)^{n-\frac{\ell}{2}}-2}{k-2}] & \text{when } \ell \textnormal{ is even},\\[15pt]
		k(k-1)^{\frac{\ell-1}{2}}[\frac{(k-1)^{n}-(k-1)^\frac{\ell-1}{2}}{k-2}] & \text{when } \ell \textnormal{ is odd.}
		\end{array} \right.
	\end{equation*}
\end{cor}

\begin{cor}\label{cor2}
For natural number $n$ and $k$ such that $k \geq 3$, we have that
\begin{equation*}
\binom{1 + \frac{k[(k - 1)^{n} - 1]}{k - 2}}{2}= \sum^{n-1}_{\ell=0}k(k-1)^{\ell}[\frac{(k-1)^{n}-(k-1)^{\ell}}{k-2}]+\sum^{n}_{\ell=1}\frac{k(k-1)^{2\ell-1}}{2} [\frac{k(k-1)^{n-l}-2}{k-2}]
\end{equation*}
\end{cor}

\subsection{Wiener Index and Average Distance}\label{ap1}
In this subsection,We have linked our main problem to distance in Cayaley Tree Dendrimer. Using the results obtained in Theorem \ref{mainthm},  Corollary \ref{cor1} and \ref{cor2}, we have found the Wiener index and average distance of $T_{n,k}$. We obtain Corollaries \ref{weiofden} and \ref{muofden}. However, we may need Theorem \ref{totaldistanceoftree} and the proof of this theorem is given in Section \ref{proofs}.
\vskip 5 pt

\begin{thm}\label{totaldistanceoftree}
Let $T$ be a tree having the diameter $diam(T)$. Then
$$\displaystyle{\sum_{\{u, v\} \subseteq V(T)}d_{T}(u, v)} = \displaystyle{\sum^{diam(T)}_{\ell = 1}\ell n_{\ell}(T)}.$$
\end{thm}
\vskip 5 pt

\noindent By Corollary \ref{cor1}, we have that
\vskip 5 pt

\begin{align}\label{eq}
\displaystyle{\sum^{2n}_{\ell = 1}\ell n_{\ell}(T_{n, k})}= &\sum^{n-1}_{l=0}(2l+1)k(k-1)^{l}[\frac{(k-1)^{n}-(k-1)^{l}}{k-2}] +\notag\\
&\sum^{n}_{l=1}(2l)\frac{k(k-1)^{2l-1}}{2} [\frac{k(k-1)^{n-l}-2}{k-2}].
\end{align}

\noindent As $diam(T_{n, k}) = 2n$, by (\ref{eq}) and Theorem \ref{totaldistanceoftree}, we immediately obtain the following corollaries.
\vskip 5 pt

\begin{cor}\label{weiofden}
Let $T_{n, k}$ be the dendrimer with the Weiner index $W(T_{n, k})$. Then
$$W(T_{n, k}) = \sum^{n-1}_{l=0}(2l+1)k(k-1)^{l}[\frac{(k-1)^{n}-(k-1)^{l}}{k-2}]+\sum^{n}_{l=1}(2l)\frac{k(k-1)^{2l-1}}{2} [\frac{k(k-1)^{n-l}-2}{k-2}].$$
\end{cor}
\vskip 5 pt

\noindent It is worth noting that the right hand side of the equation in Corollary \ref{weiofden} can be simplified to Equation (\ref{weiner1}).
\vskip 5 pt

\begin{cor}\label{muofden}
Let $T_{n, k}$ be the dendrimer with the average distance $\mu(T_{n, k})$. Then
$$\mu(T_{n, k}) = \frac{\sum^{n-1}_{l=0}(2l+1)k(k-1)^{l}[\frac{(k-1)^{n}-(k-1)^{l}}{k-2}]+\sum^{n}_{l=1}(2l)\frac{k(k-1)^{2l-1}}{2} [\frac{k(k-1)^{n-l}-2}{k-2}]}{\displaystyle{\binom{1 + \frac{k[(k - 1)^{n} - 1]}{k - 2}}{2}}}.$$
\end{cor}

\subsection{Medium Domination}\label{ap2}
Motivated by \cite{GS,GV3,DP}, we generalize their results to  $\varsigma$-medium domination of $T_{n, k}$. For a graph $G$ of order $n$ and for some $2 \leq \varsigma \leq diam(G)$, the $\varsigma$\emph{-medium domination number} $\gamma_{\varsigma MD}(G)$ of $G$ is defined as 
$$\gamma_{\varsigma MD}(G) = \frac{\varsigma(G)}{{n\choose 2}}$$ 
\noindent where 
$$\varsigma(G)=\sum^{\varsigma}_{\ell = 1}n_{\ell}(G) ,$$
\noindent the sum of all paths whose lengths less than or equal to $\varsigma$. Hence, when $G$ is a dendrimer $T_{n, k}$, we obtain the $\varsigma$-medium domination number of $T_{n, k}$ as follow:

\begin{cor}
Let $T_{n, k}$ be the dendrimer with the $\varsigma$-medium domination number $\gamma_{\varsigma MD}(G)$. Then
\end{cor}
$$\gamma_{\varsigma MD}(T_{n, k})= \frac{\varsigma(T_{n, k})}{{|V(T_{n, k})|\choose 2}}$$ 
\noindent where
\begin{align*}
\varsigma(T_{n, k}) =\sum^{s}_{\ell=0}k(k-1)^{\ell}[\frac{(k-1)^{n}-(k-1)^{\ell}}{k-2}]+\sum^{\lfloor \frac{\varsigma}{2}\rfloor}_{\ell=1}\frac{k(k-1)^{2\ell-1}}{2} [\frac{k(k-1)^{n-\ell}-2}{k-2}]
\end{align*}
\noindent and
\begin{equation*}
	s = \left\{
	\begin{array}{rl}
		\lfloor \frac{\varsigma}{2}\rfloor & \text{when } \varsigma \textnormal{ is odd},\\[15pt]
		\lfloor \frac{\varsigma}{2}\rfloor - 1 & \text{when } \varsigma \textnormal{ is even.}
		\end{array} \right.
	\end{equation*}


\section{Preliminaries}
In this section, we provide some results hat are used in establishing our main theorems. We begin with a simple but yet useful formula for geometric series. For a geometric series $S_{n} = a + ar + ar^{2} + \cdots + ar^{n - 1}$, we have that 

$$S_{n} = \sum^{n - 1}_{i = 0}ar^{i} = \frac{a(1-r^n)}{1-r}$$

\noindent where $n$ is the number of terms, $a$ is the coefficient and $r \neq 1$ is the common ratio.

\vskip 5 pt

\noindent Further, for $T_{n, k}$, we may have the following formulas by simple counting arguments and geometric series,
\vskip 5 pt

\indent the total number of vertices of degree $k$ are  $\frac{k(k-1)^{n-1}-2}{k-2}$,
\vskip 5 pt

\indent the total number of vertices of degree 1 (i.e. pendent vertices) is equal to $k(k-1)^{n-1}$,
\vskip 5 pt

\indent the total number of vertices is equal to $1+\frac{k[(k-1)^{n}-1]}{k-2}$,
\vskip 5 pt

\noindent and
\vskip 5 pt

\indent the total number of edges is equal to $\frac{k[(k-1)^{n}-1]}{k-2}$.
\vskip 5 pt

\section{Proofs}\label{proofs}
In this section, we give the proofs of Theorem \ref{mainthm}, Corollary \ref{cor1}, Corollary \ref{cor2} and Theorem \ref{totaldistanceoftree}.

\subsection{Proof of Theorem \ref{mainthm}}
\indent To prove this theorem, we need to establish Lemmas \ref{lem1} and \ref{lem2} which are the exact formulas of $n^{1}_{\ell}(T_{n, k})$ and $n^{2}_{\ell}(T_{n, k})$.
\vskip 5 pt

\begin{lem}\label{lem1}
for $n, k \geq 1$ and $1 \leq \ell \leq 2n$, we let $n^{1}_{\ell}(T_{n, k})$ be the number of paths of length $\ell$ of $T_{n, k}$ having exactly one end vertex as a leaf of $T_{n, k}$. Then
    \begin{equation*}
	n^{1}_{\ell}(T_{n, k}) = \left\{
	\begin{array}{rl}
		k(k - 1)^{n + \frac{\ell}{2} - 2} & \text{when } \ell \textnormal{ is even},\\[15pt]
		k(k - 1)^{n + \frac{\ell - 1}{2} - 1} & \text{when } \ell \textnormal{ is odd.}
		\end{array} \right.
	\end{equation*}
\end{lem}

\proof First, we let $r$ be the root and let $x$ be an arbitrary leaf of the graph $T_{n, k}$. Further, for $0 \leq j \leq n$, we let 
\vskip 5 pt

\indent $L_{j}:$ the set of all vertices of $T_{n, k}$ at distance $i$ from $r$,
\vskip 5 pt

\noindent and
\vskip 5 pt

\indent $\mathcal{P}_{x}:$ the family of all paths of $T_{n, k}$ starting from $x$ and the other end vertex is not\\
\indent $~~~~$ a leaf of $T_{n, k}$.
\vskip 5 pt

\noindent We distinguish two cases according to the value of $\ell$.
\vskip 5 pt

\noindent \textbf{Case 1:} $1 \leq\ell \leq n$.

\noindent For a path $P \in \mathcal{P}_{x}$, we let
\vskip 5 pt

\indent $min(P) = \min\{j : V(P) \cap L_{j} \neq \emptyset\}$.
\vskip 5 pt

\noindent Further, for $0 \leq i \leq \lfloor\frac{\ell - 1}{2}\rfloor$, we let
\vskip 5 pt

\indent $\mathcal{P}_{x, i} = \{P \in \mathcal{P}_{x} : min(P) = n - \ell + i\}$.
\vskip 5 pt

\noindent It can be observed that $\mathcal{P}_{x, 0}, \mathcal{P}_{x, 1}, ..., \mathcal{P}_{x, \lfloor\frac{\ell - 1}{2}\rfloor}$ partition $\mathcal{P}_{x}$.
\vskip 5 pt

\indent When $i = 0$, we have that $|\mathcal{P}_{x, 0}| = 1$ as there is exactly one path of length $\ell$ starting from $x$, goes through vertices in $L_{n - 1}, L_{n - 2}, ..., L_{n - \ell + 1}$ and terminates in $L_{n - \ell}$.
\vskip 5 pt

\indent For each $1 \leq i \leq \lfloor\frac{\ell - 1}{2}\rfloor$, all the paths in $\mathcal{P}_{x, i}$ start from $x$ and go trough vertices in $L_{n - 1}, ..., L_{n - \ell + i + 1}, L_{n - \ell + i}$ with exactly one possibility. We may let $y \in L_{n - \ell + i + 1}$ and $z \in L_{n - \ell + i}$ be the vertices that are in all the paths. Then, from the vertex $z$, all the paths move back to $L_{n - \ell + i + 1}, ..., L_{n - \ell + 2i}$. As $y$ is already in every of such path, there are $k - 2$ possibilities for all the paths in $P_{x, i}$. Further, there are $k - 1$ possibilities for all the paths to pass each of $L_{n - \ell + i + 2}, ..., L_{n - \ell + 2i}$. Hence, 
\begin{align*}
    |\mathcal{P}_{x, i}| = (k - 2)(k - 1)^{i - 1}
\end{align*}
\noindent which implies that
\begin{align*}
    |\mathcal{P}_{x}| &= |\mathcal{P}_{x, 0}| + |\mathcal{P}_{x, 1}| + \cdots + |\mathcal{P}_{x, \lfloor\frac{\ell - 1}{2}\rfloor}|\\
                      &= 1 + (k-2) + (k-2)(k-1) + \cdots + (k-2)(k-1)^{\lfloor\frac{\ell-3}{2}\rfloor}.
\end{align*}

\noindent After simplifying this geometric series, we get 
\begin{equation*}  
|\mathcal{P}_{x}| =(k-1)^{\lfloor\frac{\ell-1}{2}\rfloor}
\end{equation*}
\noindent and this proves Case 1.
\vskip 5 pt


\noindent \textbf{Case 2:} $\ell = n + 1\leq l \leq 2n$.

\noindent In this case, we let
\vskip 5 pt

\indent $\mathcal{R}_{x} = \{P \in \mathcal{P}_{x} : r \in V(P)\}$
\vskip 5 pt

\noindent and
\vskip 5 pt

\indent $\mathcal{S}_{x} = \{P \in \mathcal{P}_{x} : r \notin V(P)\}$.
\vskip 5 pt

\indent We first count the number of paths in $\mathcal{R}_{x}$. All the paths in $\mathcal{R}_{x}$ start from $x$ and pass to the root $r$ with one possibilities. Then, from $r$, all the paths pass trough  $L_{1}, ..., L_{\ell - n - 1}$ and terminate in $L_{\ell - n}$, each of which with the possibilities $k - 1$. Thus, $|\mathcal{R}_{x}| = (k - 1)^{\ell - n}$.
\vskip 5 pt

\indent Next, we count the number of paths in $\mathcal{S}_{x}$ by similar arguments as in Case 1. For $\ell - n + 1 \leq i \leq \lfloor\frac{\ell - 1}{2}\rfloor$, we let
\vskip 5 pt

\indent $\mathcal{S}_{x, i} = \{P \in \mathcal{S}_{x} : min(P) = n - \ell + i\}$.
\vskip 5 pt

\noindent Clearly, $\mathcal{S}_{x, \ell - n + 1}, ..., \mathcal{S}_{x, \lfloor\frac{\ell - 1}{2}\rfloor}$ partitions $\mathcal{S}_{x}$.
\vskip 5 pt

\indent For each $\ell - n + 1 \leq i \leq \lfloor\frac{\ell - 1}{2}\rfloor$, all paths in $\mathcal{S}_{x, i}$ start from $x$ pass trough $L_{n - 1}, ..., L_{n - \ell + i + 1}$ to $L_{n - \ell + i}$ with one possibility. Then, the paths pass back to $L_{n - \ell + i + 1}$ with $k - 2$ possibilities and continue in $L_{n - \ell + i + 2}$ until terminating in $L_{n - \ell + 2i}$, each of which with $k - 1$ possibilities. Thus 
\begin{align*}
    |\mathcal{S}_{x, i}| = (k - 2)(k - 1)^{i - 1}
\end{align*}
\noindent which implies that
\begin{align*}
    |S_{x}| &= |\mathcal{S}_{x, \ell - n + 1}| +  \cdots + |\mathcal{S}_{x, \lfloor\frac{\ell - 1}{2}\rfloor}|\\
            &= (k - 2)(k - 1)^{\ell - n} + (k - 2)(k - 1)^{\ell - n + 1} +\cdots + (k - 2)(k - 1)^{\lfloor\frac{\ell - 3}{2}\rfloor}\\
            &= (k - 2)\Big(\frac{(k - 1)^{\lfloor\frac{\ell - 1}{2}\rfloor} - 1}{k - 2} - \frac{(k - 1)^{\ell - n} - 1}{k - 2}\Big)\\
            &= (k - 1)^{\lfloor\frac{\ell - 1}{2}\rfloor} - (k - 1)^{\ell - n}.
\end{align*}

\noindent Hence,
\begin{align*}
    |\mathcal{P}_{x}| = |\mathcal{R}_{x}| + |\mathcal{S}_{x}| = (k - 1)^{\lfloor\frac{\ell - 1}{2}\rfloor}
\end{align*}
\noindent and this proves Case 2.
\vskip 5 pt

\indent In both cases, we have that $|\mathcal{P}_{x}| =(k-1)^{\lfloor\frac{\ell-1}{2}\rfloor}$. As $x$ is an arbitrary leaf of $T_{n, k}$ and $T_{n, k}$ has $k(k - 1)^{n - 1}$ leaves, it follows that
\begin{equation*}
	n^{1}_{\ell}(T_{n, k}) = \left\{
	\begin{array}{rl}
		k(k - 1)^{n + \frac{\ell}{2} - 2} & \text{when } \ell \textnormal{ is even},\\[15pt]
		k(k - 1)^{n + \frac{\ell - 1}{2} - 1} & \text{when } \ell \textnormal{ is odd}
		\end{array} \right.
	\end{equation*}
\noindent and this proves Lemma \ref{lem1}.
\qed


\begin{lem}\label{lem2}
Let $n^{2}_{\ell}(T_{n, k})$ be the number of paths of length $\ell$ that starts and end on a leaf vertex of the graph $T_{n,k}$. Then, for $n,k \geq1$,
    \begin{equation*}
	n^{2}_{\ell}(T_{n, k}) = \left\{
	\begin{array}{rl}
		k(k - 1)^{n+\frac{\ell}{2}- 3}\binom{k - 1}{2} & \text{when } 2 \leq \ell \leq 2n-2\\[15pt]
		(k - 1)^{\ell - 2}\binom{k}{2} & \text{when } \ell = 2n.
		\end{array} \right.
	\end{equation*}
\end{lem}

\proof First, we let
\vskip 5 pt

\indent $\mathcal{Q}_{\ell}:$ the family of paths of length $\ell$ of $T_{n, k}$ whose both end vertices are leaves of $T_{n, k}$.
\vskip 5 pt

\noindent Clearly, $\ell$ must be even. For a path $P \in \mathcal{Q}_{\ell}$, we let $x_{P}$ be the center of $P$ which the distance from $x_{p}$ to the end vertices of $P$ are both equal to $\frac{\ell}{2}$. We distinguish $2$ cases according to the value of $\ell$.
\vskip 5 pt

\noindent \textbf{Case 1:}  $2 \leq \ell \leq 2n-2$.

\indent It can be observed that every path in $\mathcal{Q}_{\ell}$ has the center in $L_{n - \frac{\ell}{2}}$. Let $x$ be a vertex in $L_{n - \frac{\ell}{2}}$. There are $k - 1$ neighbors of $x$ in $L_{n - \frac{\ell}{2} + 1}$. Each pair of these $k - 1$ neighbors can be passed by a path in $\mathcal{Q}_{\ell}$. Hence, there are $\binom{k - 1}{2}$ possibilities for the paths in $\mathcal{Q}_{\ell}$. We may let $x_{1}$ and $x_{2}$ be a pair among these $\binom{k - 1}{2}$ possibilities. There are $(k - 1)^{\frac{\ell}{2} - 1}$ paths from each of $x_{1}$ and $x_{2}$ to the leaves of $T_{n, k}$. Hence, there are 
\begin{equation*}
\binom{k - 1}{2}(k - 1)^{\frac{\ell}{2} - 1}(k - 1)^{\frac{\ell}{2} - 1} = (k - 1)^{\ell - 2}\binom{k-1}{2}
\end{equation*}
\noindent paths whose center is $x$ and both end vertices are leaves. Since $x$ is arbitraty and there are $k(k - 1)^{n - \frac{\ell}{2} - 1}$ vertices in $L_{n - \frac{\ell}{2}}$, it follows that
\begin{align*}
    n^{2}_{\ell}(T_{n, k}) = |\mathcal{Q}_{\ell}| = k(k - 1)^{n - \frac{\ell}{2} - 3}\binom{k-1}{2}.
\end{align*}

\vskip 5 pt

\noindent \textbf{Case 2:} $\ell = 2n$

\indent In this case, the root $r$ is the center of all paths in $\mathcal{Q}_{2n}$. There are $\binom{k}{2}$ possibilities for the paths in $\mathcal{Q}_{2n}$ to pass these vertices. Similarly, we let $x_{1}$ and $x_{2}$ be a pair among these $\binom{k}{2}$ possibilities. There are $(k - 1)^{\frac{\ell}{2} - 1}$ paths from each of $x_{1}$ and $x_{2}$ to the leaves of $T_{n, k}$. Hence, 
\begin{equation*}
n^{2}_{\ell}(T_{n, k}) = |\mathcal{Q}_{2n}| = (k - 1)^{\ell - 2}\binom{k}{2}
\end{equation*}
\noindent and this proves Lemma \ref{lem2}.
\qed
\vskip 5 pt

\indent Now we are ready to prove Theorem \ref{mainthm}.
\vskip 5 pt

\noindent \textbf{Proof of Theorem \ref{mainthm}} Recall that the graph $T_{n, k}$ can be constructed from $T_{n - 1, k}$ by introducing $k - 1$ vertices to each leaf, and joining these $k - 1$ vertices to the leaf.
We have considered two cases.
\vskip 5 pt

\noindent \textbf{Case 1:} $\ell$ is an even number.

\indent Every path of length $\ell$ in this case is either (\emph{i}) lies completely in $T_{n-1, k}$, (\emph{ii}) can be formed from a path of length $\ell - 1$ whose exactly one end vertex is a leaf of $T_{n - 1, k}$ or (\emph{iii}) can be formed from a path of length $\ell - 2$ whose both end vertices are at the leaves of $T_{n - 1, k}$. The Case (\emph{i}) gives $n_{\ell}(T_{n - 1, k})$ paths of length $\ell$ while the Case (\emph{ii}) gives $(k - 1)n^{1}_{\ell}(T_{n - 1, k})$ paths of length $\ell$ as the end vertex at a leaf of $T_{n - 1, k}$ can be extended with $k - 1$ ways. Finally, the Case (\emph{iii}) gives $(k - 1)^{2}n^{2}_{\ell - 2}(T_{n - 1, k})$ paths as every path of length $\ell - 2$ whose both end vertices are at the leaves of $T_{n - 1, k}$ can be extended to the path of length $\ell$ by $(k - 1)^{2}$ ways, $k - 1$ for each end vertex. 
Thus, we have the following recursive formula

\begin{equation*}
n_{\ell}(T_{n, k}) = (k - 1)n^{1}_{\ell - 1}(T_{n - 1, k}) + (k - 1)^{2}n^{2}_{\ell - 2}(T_{n - 1, k}) + n_{\ell}(T_{n - 1, k}).
\end{equation*}

\noindent This proves Case 1.
\vskip 8 pt

\noindent \textbf{Case 2:} $\ell$ is an odd number

\indent Similarly, every path of length $\ell$ in this case is either (\emph{i}) lies completely in $T_{n-1, k}$, (\emph{ii}) can be formed from a path of length $\ell - 1$ whose exactly one end vertex is a leaf of $T_{n - 1, k}$ or (\emph{iii}) can be formed from a path of length $\ell - 1$ whose both end vertices are at the leaves of $T_{n - 1, k}$. The Case (\emph{i}) gives $n_{\ell}(T_{n - 1, k})$ paths while the Case (\emph{ii}) gives $(k - 1)n^{1}_{\ell}(T_{n - 1, k})$ paths. For the Case (\emph{iii}), we can only extend these paths of length $\ell - 1$ in $T_{n - 1, k}$ to be a path of length $\ell$ by extending only one end vertex, $k - 1$ ways for each end vertex. Thus there are $2 (k - 1)n^{2}_{\ell - 1}(T_{n - 1, k})$ paths in this case. Thus, we have a recursive formula

\begin{equation*}
n_{\ell}(T_{n, k}) = (k - 1)n^{1}_{\ell - 1}(T_{n - 1, k}) +2 (k - 1)n^{2}_{\ell - 1}(T_{n - 1, k}) + n_{\ell}(T_{n - 1, k}).
\end{equation*}

\noindent This proves Case 2 and completes the proof of our theorem.
\qed

\subsection{Proof of Corollary \ref{cor1}}
We distinguish two cases according to the parity of $\ell$.
\vskip 5 pt

\noindent \textbf{Case 1:} $\ell$ is an even number.

By Theorem \ref{mainthm}, we have that

\begin{align}
n_{\ell}(T_{n, k}) &= (k - 1)n^{1}_{\ell - 1}(T_{n - 1, k}) + (k - 1)^{2}n^{2}_{\ell - 2}(T_{n - 1, k}) + n_{\ell}(T_{n - 1, k})\notag\\
n_{\ell}(T_{n-1, k}) &= (k - 1)n^{1}_{\ell - 1}(T_{n - 2, k}) + (k - 1)^{2}n^{2}_{\ell - 2}(T_{n - 2, k}) + n_{\ell}(T_{n - 2, k})\notag\\
n_{\ell}(T_{n-2, k}) &= (k - 1)n^{1}_{\ell - 1}(T_{n - 3, k}) + (k - 1)^{2}n^{2}_{\ell - 2}(T_{n - 3, k}) + n_{\ell}(T_{n - 3, k})\notag\\
& \vdots \notag\\
n_{\ell}(T_{\frac{\ell}{2}+1, k}) &= (k - 1)n^{1}_{\ell - 1}(T_{\frac{\ell}{2}, k}) + (k - 1)^{2}n^{2}_{\ell - 2}(T_{\frac{\ell}{2}, k}) + n_{\ell}(T_{\frac{\ell}{2}, k})\notag\\
n_{\ell}(T_{\frac{\ell}{2}, k}) &= (k - 1)n^{1}_{\ell - 1}(T_{\frac{\ell - 2}{2}, k}) + (k - 1)^{2}n^{2}_{\ell - 2}(T_{\frac{\ell - 2}{2}, k}) + n_{\ell}(T_{\frac{\ell - 2}{2}, k}).\notag
\end{align}

\noindent As $\frac{\ell - 2}{2} < \frac{\ell}{2}$, we have $n_{\ell}(T_{\frac{\ell - 2}{2}, k}) = 0$. Further, $n^{1}_{\ell- 1}(T_{\frac{\ell - 2}{2}, k}) = 0$ because $\frac{\ell - 2}{2} < \frac{\ell - 1}{2}$. Thus, summing the above equations we have

\begin{equation}\label{rf}
   n_{\ell}(T_{n, k}) = (k - 1)\sum^{n-1}_{i=\frac{\ell}{2}}{n^{1}_{\ell - 1}(T_{i, k})} + (k - 1)^{2}\sum^{n-1}_{j=\frac{\ell-2}{2}}{n^{2}_{\ell - 2}(T_{j, k})}. 
\end{equation}

\noindent By Lemma \ref{lem1} when $\ell - 1$ is odd, we have that

\begin{align}
\sum^{n-1}_{i=\frac{\ell}{2}}{n^{1}_{\ell - 1}(T_{i, k})} &=\sum^{n-1}_{i=\frac{\ell}{2}}{k(k-1)^{i+\frac{\ell-2}{2}-1}}\notag\\
&= k(k-1)^{\frac{\ell}{2}-2}[\sum^{n-1}_{i=0}{(k-1)^i}-\sum^{\frac{\ell}{2}-1}_{i=0}{(k-1)^i}].\notag
\end{align}

\noindent By Geometric Series, we have that

\begin{equation}\label{l1}
  \sum^{n-1}_{i=\frac{\ell}{2}}{n^{1}_{\ell - 1}(T_{i, k})}= k(k-1)^{\frac{\ell}{2}-2}[\frac{(k-1)^{n}-(k-1)^{\frac{\ell}{2}}}{k-2}].  
\end{equation}

\noindent For the sum $\sum^{n-1}_{j=\frac{\ell-2}{2}}{n^{2}_{\ell - 2}(T_{j, k})}$, we may split the first term as

\begin{equation*}
    \sum^{n-1}_{j=\frac{\ell-2}{2}}{n^{2}_{\ell - 2}(T_{j, k})}=n^{2}_{\ell - 2}(T_{\frac{\ell-2}{2}, k})+\sum^{n-1}_{j=\frac{\ell}{2}}{n^{2}_{\ell - 2}(T_{j, k})}.
\end{equation*}

\noindent By Lemma \ref{lem2}, we have that

\begin{align}
\sum^{n-1}_{j=\frac{\ell-2}{2}}{n^{2}_{\ell - 2}(T_{j, k})} 
&=(k-1)^{\ell} {k \choose 2}+\sum^{n-1}_{j=\frac{\ell}{2}}{k(k-1)^{\frac{\ell}{2}-4+j}}{k-1 \choose 2}\notag\\
&=(k-1)^{\ell} {k \choose 2}+k(k-1)^{\frac{\ell}{2}-4}{k-1 \choose 2}\sum^{n-1}_{j=\frac{\ell}{2}}(k-1)^{j}.\notag
\end{align}

\noindent By Geometric Series, we have that

\begin{equation}\label{l2}
    \sum^{n-1}_{j=\frac{\ell-2}{2}}{n^{2}_{\ell - 2}(T_{j, k})}=(k-1)^{\ell-4} {k \choose 2}+k(k-1)^{\frac{\ell}{2}-4}{k-1 \choose 2}[\frac{(k-1)^{n}-(k-1)^{\frac{\ell}{2}}}{k-2}].
\end{equation}
\vskip 5 pt

\noindent Putting values from Equation (\ref{l1}) and (\ref{l2}) into Equation (\ref{rf}) and simplifying, we get 
\begin{equation}\label{final formula 1}
n_{\ell}(T_{n, k})= \frac{k(k-1)^{\ell-1}}{2} [\frac{k(k-1)^{n-\frac{\ell}{2}}-2}{k-2}].
\end{equation}

\noindent This proves Case 1.

\vskip 8 pt

\noindent \textbf{Case 2:} $\ell$ is an odd number

\noindent By Theorem \ref{mainthm}, we have that

\begin{align}
n_{\ell}(T_{n, k})    &= (k - 1)n^{1}_{\ell - 1}(T_{n - 1, k}) +2 (k - 1)n^{2}_{\ell - 1}(T_{n - 1, k}) + n_{\ell}(T_{n - 1, k})\notag\\
n_{\ell}(T_{n-1, k})  &= (k - 1)n^{1}_{\ell - 1}(T_{n - 2, k}) + 2(k - 1)n^{2}_{\ell - 1}(T_{n - 2, k}) + n_{\ell}(T_{n - 2, k})\notag\\
n_{\ell}(T_{n-2, k})  &= (k - 1)n^{1}_{\ell - 1}(T_{n - 3, k}) + 2(k - 1)n^{2}_{\ell - 1}(T_{n - 3, k}) + n_{\ell}(T_{n - 3, k})\notag\\
&\vdots \notag\\
n_{\ell}(T_{\frac{\ell+1}{2}, k}) &= (k - 1)n^{1}_{\ell- 1}(T_{\frac{\ell - 1}{2}, k}) + 2(k - 1)n^{2}_{\ell - 1}(T_{\frac{\ell - 1}{2}, k}) + n_{\ell}(T_{\frac{\ell - 1}{2}, k}).\notag
\end{align}

\noindent Since $\frac{\ell - 1}{2} < \frac{\ell}{2}$, it follows that $n_{\ell}(T_{\frac{\ell - 1}{2}, k}) = 0$. Further, $n^{1}_{\ell- 1}(T_{\frac{\ell - 1}{2}, k}) = 0$ because every path of length $\ell - 1$ always has both end vertices at leaves of $T_{\frac{\ell - 1}{2}, k}$. Thus, summing the above equations we have

\begin{equation}\label{rf2}
    n_{\ell}(T_{n, k}) = (k - 1)\sum^{n-1}_{i=\frac{\ell+1}{2}}{n^{1}_{\ell - 1}(T_{i, k})} + 2(k - 1)\sum^{n-1}_{j=\frac{\ell+1}{2}}{n^{2}_{\ell - 1}(T_{j, k})}+2(k-1){n^{2}_{\ell - 1}}(T_{\frac{\ell-1}{2}, k}).
\end{equation} 

\noindent By Lemma \ref{lem1} when $\ell - 1$ is even, we have that
\begin{align}
\sum^{n-1}_{i=\frac{\ell+1}{2}}{n^{1}_{\ell - 1}(T_{i, k})}     &=\sum^{n-1}_{i=\frac{\ell+1}{2}}{k(k-1)^{i+\frac{\ell-1}{2}-2}}\notag\\
&=k(k-1)^{\frac{\ell-1}{2}-2}\sum^{n-1}_{i=\frac{\ell+1}{2}}{(k-1)^{i}}\notag\\
&=k(k-1)^{\frac{\ell-1}{2}-2}[\sum^{n-1}_{i=0}{(k-1)^{i}}-\sum^{\frac{\ell-1}{2}}_{i=0}{(k-1)^{i}}].\notag
\end{align} 

\noindent Hence, we have by Geometric Series that
 \begin{equation}\label{l11}
     \sum^{n-1}_{i=\frac{\ell+1}{2}}{n^{1}_{\ell - 1}(T_{i, k})}=k(k-1)^{\frac{\ell-1}{2}-2}[\frac{(k-1)^{n}-(k-1)^{\frac{\ell+1}{2}}}{k-2}].
 \end{equation} 
 
\noindent Further, we have by Lemma \ref{lem2} that
\begin{align}
\sum^{n-1}_{j=\frac{\ell+1}{2}}{n^{2}_{\ell - 1}(T_{j, k})} &=\sum^{n-1}_{j=\frac{\ell+1}{2}}k(k-1)^{\frac{\ell+1}{2}-2+j}{k -1\choose 2} \notag\\
&=k(k-1)^{\frac{\ell+1}{2}-2}{k -1\choose 2}\sum^{n-1}_{j=\frac{\ell+1}{2}}(k-1)^{j}.\notag
\end{align}

\noindent We have by Geometric Series that
\begin{equation}\label{l22}
    \sum^{n-1}_{j=\frac{\ell+1}{2}}{n^{2}_{\ell - 1}(T_{j, k})}=k(k-1)^{\frac{\ell+1}{2}-2}{k -1\choose 2}[\frac{(k-1)^{n}-(k-1)^{\frac{\ell+1}{2}}}{k-2}].
\end{equation} 

\noindent Putting values from Equations (\ref{l11}) and (\ref{l22}) into Equation (\ref{rf2}) and simplifying, we get
\newline \begin{equation}\label{final formula 2}
n_{\ell}(T_{n, k}) = k(k-1)^{\frac{\ell-1}{2}}[\frac{(k-1)^{n}-(k-1)^\frac{\ell-1}{2}}{k-2}].
\end{equation}

\noindent This proves Case 2 and completes the proof of Corollary \ref{cor1}.
\qed

\subsection{Proof of Corollary \ref{cor2}}
We let $\binom{V(T_{n, k})}{2}$ be the set of all sets of two vertices of $T_{n, k}$. Namely, 
\begin{align*}
\binom{V(T_{n, k})}{2} = \{\{u, v\} : u, v \in V(T_{n, k})\} 
\end{align*}
and 
\begin{align*}
\Big|\binom{V(T_{n, k})}{2}\Big| = \binom{|V(T_{n, k})|}{2} = \binom{1 + \frac{k[(k - 1)^{n} - 1]}{k - 2}}{2}.
\end{align*}

\noindent Construct the $(0, 1)$-matrix whose rows are the pairs $\{u, v\}$ of $\binom{V(T_{n, k})}{2}$, columns are the path length $\ell$ for all $1 \leq \ell \leq 2n$ and the entries $a_{\{u, v\}, \ell}$ are defined as follows:

\begin{equation*}
	a_{\{u, v\}, \ell} = \left\{
	\begin{array}{rl}
		1 & \text{if } d_{T}(u, v) = \ell,\\[15pt]
		0 & \text{otherwise.}
		\end{array} \right.
\end{equation*}

\indent We first consider Row $\{u, v\}$. There is exactly one column, $\ell$ say, such that 
$$a_{\{u, v\}, \ell} = 1$$
\noindent but
$$a_{\{u, v\}, j} = 0$$
\noindent for all $j \in \{1, ..., 2n\} \setminus \{\ell\}$. Thus, the summation of all entries in this matrix is 
\begin{align*}
    \sum_{\{u, v\} \subseteq V(T)}1 = \Big|\binom{V(T_{n, k})}{2}\Big| = \binom{1 + \frac{k[(k - 1)^{n} - 1]}{k - 2}}{2}.
\end{align*}

\vskip 5 pt

\indent We then consider Column $\ell$. By the definition of $n_{\ell}(T_{n, k})$, there are $n_{\ell}(T_{n, k})$ rows whose entries are equal to $1$ while the entries of the other rows are all $0$. Hence, the summation of all entries of Column $\ell$ is equal to $n_{\ell}(T)$ implying that the summation of all etries in this matrix is $\sum^{2n}_{\ell = 1}n_{\ell}(T)$.
\vskip 5 pt

\indent By the counting two way principle, we have that

$$\binom{1 + \frac{k[(k - 1)^{n} - 1]}{k - 2}}{2} = \sum^{2n}_{\ell = 1}n_{\ell}(T_{n, k}).$$

\noindent This proves Corollary \ref{cor2}.

\subsection{Proof of Theorem \ref{totaldistanceoftree}}
We prove this theorem by similar argument as in the proof of Corollary \ref{cor2}. First, we let $T$ be a tree with the diameter $diam(T) = t$. We let $\binom{V(T)}{2}$ be the set of all sets of two vertices of $T$. Construct the matrix whose rows are the pairs $\{u, v\}$ of $\binom{V(T)}{2}$, columns are the path length $\ell$ for all $1 \leq \ell \leq t$ and the entries $a_{\{u, v\}, \ell}$ are defined as follows:

\begin{equation*}
	a_{\{u, v\}, \ell} = \left\{
	\begin{array}{rl}
		\ell & \text{if } d_{T}(u, v) = \ell,\\[15pt]
		0 & \text{otherwise.}
		\end{array} \right.
\end{equation*}

\noindent Thus, in the Column $j$, all the entries are either $j$ or $0$.
\vskip 5 pt

\indent We first consider Row $\{u, v\}$. There is exactly one column, $\ell$ say, such that 
$$a_{\{u, v\}, \ell} = \ell = d_{T}(u, v)$$
\noindent but
$$a_{\{u, v\}, j} = 0$$
\noindent for all $j \in \{1, ..., t\} \setminus \{\ell\}$. Thus, the summation of all entries of Row $\{u, v\}$ is equal to $\ell = d_{T}(u, v)$ implying that the summation of all entries in this matrix is $\sum_{\{u, v\} \subseteq V(T)}d_{T}(u, v)$.
\vskip 5 pt

\indent We then consider Column $\ell$. By the definition of $n_{\ell}(T)$, there are $n_{\ell}(T)$ rows whose entries are equal to $\ell$ while the entries of the other rows are all $0$. Hence, the summation of all entries of Column $\ell$ is equal to $\ell n_{\ell}(T)$ implying that the summation of all etries in this matrix is $\sum^{diam(T)}_{\ell = 1}\ell n_{\ell}(T)$.
\vskip 5 pt

\indent By the counting two way principle, we have that

$$\sum_{\{u, v\} \subseteq V(T)}d_{T}(u, v) = \sum^{diam(T)}_{\ell = 1}\ell n_{\ell}(T).$$

\noindent This proves Theorem \ref{totaldistanceoftree}.



\end{document}